\documentclass{amsart}

\usepackage{amsmath,amssymb}

\title[John von Neumann]{Three questions about John von Neumann}

\author{Melvyn B. Nathanson}

\address{Department of Mathematics\\Lehman College (CUNY)\\Bronx, NY 10468} 

\email{melvyn.nathanson@lehman.cuny.edu}

\date{\today}

\begin{document}

\maketitle

John von Neumann was one of the greatest scientists of the 20th century 
and one of the most complicated.  
His research career began in mathematics and mathematical physics. 
Before, during, and after World War II he was deeply involved in military work 
and was a central figure in the design and construction of atomic 
and hydrogen bombs.  
He was one of the creators of game theory and its use in economics  
and he was (perhaps the most significant of his achievements) the  
dominant figure in 
the early development of computer hardware and software.  
He was also an influential US government policy maker on defense matters and a 
strong proponent of building intercontinental ballistic missiles   and designing the 
small powerful nuclear weapons that would become the warheads of ICBMs.  

Ananyo  Bhattacharya's book,  
\emph{The Man from the Future: The Visionary Life of John von Neumann}~\cite {bhat21}, 
is in part a biography of the man, but it  is more a magisterial survey 
of von Neumann's work, a series of essays that follow, in roughly chronological order,  
the different fields in which von Neumann worked.    
The first chapters are on  mathematical logic 
 and  the foundations of quantum mechanics,  
then chapters on his work on nuclear weapons, on 
computer design, 
and on mathematical economics. 
Finally, there are chapters on his consulting for the US government and 
for the RAND corporation and on his work on cellular automata. 
This is a book worth reading.

Studying the life of von Neumann, one is led to ask three related questions: 

\begin{enumerate}
\item
Von Neumann was a great scientist.  Readers of the \emph{Notices of the American 
Mathematical Society} might ask: Was he a great mathematician?

\item
Von Neumann not only strongly and successfully lobbied for the construction 
of horrendous weapons of mass destruction and their efficient delivery on presumptive enemies,  
but he also advocated a pre-emptive nuclear strike on the Soviet Union before that country 
had the ability to retaliate.  Afterwards, he forcefully promoted the doctrine 
of Mutual Assured Destruction.    Was von Neumann evil?

\item
Von Neumann was one of the first permanent professors at the Institute for Advanced Study (IAS), 
viewed by many mathematicians as mathematical heaven, 
but, as he moved away from ``pure mathematics,'' he felt increasingly unwelcome 
and uncomfortable there.  He resigned in 1956.
Why did he quit? 
\end{enumerate}

\section{Was von Neumann a great mathematician?}

John von  Neumann was born Neumann J\' anos Lajos in Budapest on December 28, 1903.
He attended an excellent high school in Budapest, but in math and science 
he far outpaced his fellow students and teachers.   His father hired private tutors, 
first G\' abor Szeg\" o and then Michael Fekete,  to teach him advanced mathematics.  
Von Neumann's first paper, published in 1922, was a joint paper with Fekete 
on the zeros of polynomials.  In 1923 and 1925 he published important papers in set theory.  
In 1926 he received a Ph.D. in mathematics from the University of Budapest 
and also a degree in chemical engineering from ETH Zurich. 

Between 1922 and 1932, von Neumann published 45 papers on mathematics 
and mathematical physics. 
He wrote foundational papers on set theory and logic, ergodic theory, game theory, 
and  functional analysis.  These were the topics \emph{du jour}.  
He also wrote papers on physics, several with his Hungarian high school friend Eugene Wigner.  
Perhaps his most significant work in this decade was  
\emph{Mathematische Grundlagen der Quantenmechanik}, published in Berlin by Springer in 1932 
and reprinted (in German) in the United States by Dover Publications in 1943.    
Freeman Dyson~\cite{dyso13} wrote,
\begin{quotation}
Johnny's book was the first exposition of quantum mechanics that made the theory mathematically 
respectable\ldots.  [But] I was surprised to discover that nobody in the physics journals ever referred 
to Johnny's book.  So far as the physicists were concerned, Johnny did not exist\ldots.  
 \end{quotation}

Stan Ulam~\cite[p. 8]{ulam58}, von Neumann's closest friend in America, wrote: 
\begin{quotation}
Mathematicians, at the outset of their creative work, are often confronted 
by two conflicting motivations: the first is to contribute to the edifice of existing work -- it is 
there that one can be sure of gaining recognition quickly by solving outstanding problems -- 
the second is the desire to blaze new trails and to create new syntheses.  
This latter course is a more risky undertaking, the final judgment of value or success 
appearing only in the future.  In his early work, Johnny chose the first of these alternatives.  
\end{quotation}
Johnny had wanted to become famous quickly and he succeeded, 
but he had doubts about the importance of his work.  
He wrote influential papers on mathematical logic, 
but in 1931 G\" odel (citing a 1927 paper of von Neumann) proved the great 
incompleteness theorem.  Von Neumann proved the mean ergodic theorem in 1931 (published in 1932), 
but Birkhoff (aware of von Neumann's result) 
proved and published the more fundamental individual ergodic theorem in 1931.  

As World War II approached, von Neumann continued to write papers in mathematics, 
including a series of papers ``On rings of operators'' that fills the entire third volume 
of his \emph{Collected Works}~\cite{vonn61}, but he became increasingly involved 
with military-related research in many areas, including ballistics 
and the development of nuclear weapons.  He was a consultant to many 
defense installations and an important figure on the Manhattan Project in Los Alamos. 
The weapons work required extensive numerical calculations, and he soon 
became the key contributor to the theoretical framework and actual building  
of the first digital computers.  He participated in the development of the ENIAC 
at the University of Pennsylvania and led the design and construction of 
the MANIAC computer at IAS.  

In a frequently quoted passage from the essay ``The mathematician,'' 
reprinted in the first volume of his \emph{Collected Works} ~\cite{vonn61}, 
von Neumann wrote 
\begin{quotation}
As a mathematical discipline travels far from its empirical source, or still more, if it is a 
second and third generation only indirectly inspired by ideas coming from ``reality'' 
it is beset with very grave dangers.  It becomes more and more purely aestheticizing, 
more and more purely l'art pour l'art\ldots.  
[W]henever this stage is reached, 
the only remedy seems to me to be the rejuvenating return to the source: 
the re-injection of more or less directly empirical ideas.
\end{quotation}
This describes a path that mathematicians might choose, 
and the path that von Neumann chose.  As Bhattacharya notes, 
\begin{quotation}
Von Neumann constantly sought new practical fields to which he could apply his mathematical genius, 
and he seemed to choose each one with an unerring sense of its potential to revolutionize human affairs. 
\end{quotation} 
Bhattacharya follows this with Freeman Dyson's summary evaluation of von Neumann's career:
\begin{quotation}
As he moved from pure mathematics to physics to economics to engineering, 
he became steadily less deep and steadily more important.
\end{quotation}

Mark Kac~\cite[p. 125]{kac85} believed that von Neumann  was ``one of the greatest mathematicians 
of our century.''
In the academic world we have not only grade inflation but also adjective inflation.  
\emph{De mortuis nil nisi bonum}.  Obituaries often call distinguished dead professors 
``great.''  
There is, of course, a distinction between doing great work in mathematics 
and being ``one of the greatest'' mathematicians.  
In the  20th century in ``pure mathematics'' 
 there were Hilbert, Poincar\' e, Weyl, G\" odel, Banach, Weil, Erd\H os, Gel'fand, Grothendieck, \ldots. 
 Was  von Neumann in that class?

\section{Was von Neumann evil?}
During World War II, von Neumann worked enthusiastically on the design and construction of atomic bombs  
and after the war he strongly encouraged the development of hydrogen bombs.  
During the war, the fear that Germany might try to manufacture nuclear weapons 
led many American physicists and engineers to work on the Manhattan Project, though at least one 
major physicist quit when it became clear that the Germans could not produce a bomb. 
The pursuit of fusion weapons was less widely supported.  On  October 30, 1949,
 the General Advisory Committee of the Atomic Energy Commission wrote: 
\begin{quotation}
[The hydrogen bomb] is not a weapon which can be used exclusively for the destruction 
of material installations of military or semi-military purposes.  
Its use therefore carries much further than the atomic bomb itself the policy 
of extermination of civilian populations.  We all hope that by one means or another, 
the development of these weapons can be avoided\ldots.  Its use would involve 
a decision to slaughter a vast number of civilians.
\end{quotation}
Von Neumann was a super-hawk.  
Bhattacharya writes (p. 208), 
\begin{quotation}
By 1946, von Neumann was predicting that devastating nuclear war was imminent\ldots.  
His answer was preventive war - a surprise attack that would wipe out the Soviet Union's nuclear arsenal 
(and a good number of its people too) before the country was able to retaliate.  
`If you say why not bomb them tomorrow, I say why not today?' he reportedly said 
in 1950.  `If you say today at 5 o'clock, I say why not one o'clock?'\footnote{Von Neumann's 
quote is from his obituary in  Life magazine~\cite{blai57}.}

\end{quotation} 
Von Neumann abandoned the effort to launch a pre-emptive strike when it became obvious 
that the Soviet Union had enough nuclear weapons to retaliate. 

Was von Neumann wrong to have worked on the atomic bomb or on the hydrogen bomb?  
Should he be criticized for striving to belong to the scientific and policy-making elite  
that developed intercontinental ballistic missiles and American nuclear policy?
Was he wrong to have forsaken ``pure mathematics'' and mathematic physics?  

There is the dual question.  Are we wrong to pursue pure mathematics 
and not to apply at least some of our intellectual skills to help solve problems plaguing this planet?  
There is no shortage of ``real world'' problems to try to solve.
Is the mantra that even the purest mathematics 
is eventually useful only self-serving self-justification? 

The answer, of course, is that there is more than one ethical choice.   
It is easy to condemn von Neumann for some of the work he did, 
but he chose honorably to think about and work on the hardest  
political and military policy problems of his time.  
We should be grateful that we are free to choose the problems we want to think about.

\section{Why did von Neumann quit the Institute?}

Von Neumann wanted to be rich.  Money was  important to him. 
He grew up in an extremely wealthy family in Budapest, and, 
according to his friends and to his biographers, 
he was most comfortable with people who were also from wealthy families.   
He enjoyed the company of generals and admirals, senators and cabinet secretaries, 
corporate bigwigs.  

He and Eugene Wigner were young, poorly paid academics in Germany 
in 1931 when they received job offers from Oswald Veblen, who was then recruiting 
for the Princeton mathematics department. 
Wigner wrote,
\begin{quotation}
One day I received a cable offering a visiting professorship at about eight times the salary 
which I had at the Institute of Technology in Berlin.  I thought this was an error in transmission.  
John von Neumann received the same cable, so we decided maybe it was true, and we accepted. 
The pay was \$3,000 for the semester, with \$1,000 for travel--a small fortune at the time. 
\end{quotation}

In 1932, the Institute for Advanced Study announced the creation of its first school, 
the School of Mathematics.  Its first appointment was Oswald Veblen.  
The following year, the Institute 
hired James Wadell Alexander, Albert Einstein, von Neumann, and Hermann Weyl.  
In 1935, Marston Morse joined the faculty.  
The Minutes of the IAS, dated January 28, 1933, record:

\begin{quotation}
RESOLVED, That Professor John von Neumann be appointed as a Professor in the School of Mathematics 
upon the following terms: that his appointment as Professor in the School of Mathematics 
date from April 1, 1933; 
that his salary be fixed at Ten thousand Dollars (\$10,000.00) per annum\ldots. 
\end{quotation}

During World War II, many mathematicians and physicists at the Institute 
did war work for the government, 
but it was expected that after the war they would resume their ``pure'' scientific research.  
Von Neumann not only did not forsake applied  science, but, when he returned to Princeton,  
began the Electronic Computer Project (ECP) at 
the Institute, a government-funded  project with a large staff of engineers and hardware 
specialists working to construct the next generation (at the time, actually, the prototypical) computer.  
His wife Kl\' ari  von Neumann wrote, 
\begin{quotation}
There he clearly stunned, or even horrified, some of his mathematical colleagues 
of the most erudite abstraction, by openly professing his great interest in other mathematical 
tools than the blackboard and chalk or pencil and paper. His proposal to build an electronic 
computing machine under the sacred dome of the Institute, was not received with applause, 
to say the least. 
\end{quotation}
Much of the programming on the computer was devoted to solving problems related 
to the development of  hydrogen bombs.  According to Theodore Taylor, a Los Alamos bomb designer, 
``The objective was pretty specifically to be able to do the coupled hydrodynamics and radiation flow necessary for H-bombs.'' 

Most of the mathematics faculty at IAS were opposed to both the ECP and to weapons research.
George Dyson~\cite[p. 82]{dyso12} writes,  
\begin{quotation}
At the Institute for Advanced Study in early 1946, even applied mathematics was out of bounds.  
Mathematicians who had worked on applications during the war were expected to leave them behind.  
Von Neumann, however, was hooked. 
\end{quotation} 
According to Freeman Dyson,
\begin{quotation}
The mathematicians certainly knew there was classified work going on.  
They may not have known it was hydrogen bombs, but it was pretty obvious. 
And they were certainly opposed to that. 
\end{quotation}
The humanities faculty at IAS also opposed the ECP.  
``In spirit we mathematicians at the Institute would cast our lot in with the humanists,''
said Marston Morse.

After the war, von Neumann was more and more involved with high level consulting with the government, and  
was soon  living in Washington and serving as a Commissioner on  the Atomic Energy Commission.  
He seems to have felt more and more unwelcome in Princeton.  
According to his wife Kl\' ari:
\begin{quotation}
The lines were drawn and after the first flurry of excitement it became clear that we 
did not belong in Princeton any more.\ldots.  
[On May 16, 1954] while in Los Angeles on air force strategic missile business\ldots, 
von Neumann met with Paul A. Dodd, dean of letters and sciences at UCLA, who offered him 
a special interdisciplinary position, with no teaching responsibilities, as professor at large. 
\end{quotation}
The UCLA job would give von Neumann unlimited opportunities to consult with IBM and other 
companies and replenish his bank account, which had been depleted 
during his years of government service.

On August 2, 1955, von Neumann was diagnosed with advanced, metastasizing cancer.
He entered Walter Reed Hospital in Washington, DC, in March 1956.
On March 19, 1956,  he wrote to  Oppenheimer: 
\begin{quotation}
After thinking matters over, I have decided to proceed in a direction which will lead me 
away from the Institute.  In order to make matters precise and definite, I plan to resign 
from the Institute's faculty, probably as of December 31.   I will write you on this subject 
in due time, but I thought that I ought to inform you of the nature of my plans without delay.  

This does not mean that I am not considering the 22 years I spent at the Institute 
the most fruitful ones of my life and scientific career.  The association with you and my 
other colleagues at the Institute has been a great experience.  

Please consider this matter strictly confidential until we communicate further.  
\end{quotation}

Bhattacharya writes, 
\begin{quotation}
Cancer had come at a particularly cruel time.   The truth was that von Neumann had been 
unhappy at the IAS for several years before his death.  `Von Neumann, when I was there 
at Princeton, was under extreme pressure' says Beno\^ it Mandelbrot, \ldots `from mathematicians, 
who were despising him for no longer being a mathematician; by the physicists, who were despising 
him for never having been a real physicist; and by everyone for having brought to Princeton 
this collection of low-class individuals called ``programmers''\ldots.  
Von Neumann was simply being shunned.  And he was not a man to take it.' 
\end{quotation}

Von Neumann never got to California.  
He  died in Walter Reed Hospital   on February 8, 1957.

Bhattacharya's book is the best introduction to the  life and work of von Neumann. 
 
\section{Other reading} 
The literature on von Neumann is enormous.  
There is a fine biography by 
Norman Macrae, \emph{John von Neumann}~\cite{macr92}, and 
also books by his brother Nicholas A. Vonneumann~\cite{vonn87} 
and  his daughter Marina von Neumann Whitman~\cite{whit12}.  
Ulam's autobiography~\cite{ulam76} contains many stories about von Neumann. 

There are serious studies on almost every aspect of von Neumann's work.  
For the history of computers, 
there are books by Goldstine~\cite{gold72} and Dyson~\cite{dyso12}.  
On the history of the development of nuclear weapons, see Rhodes~\cite{rhod87,rhod95}. 
An issue of the \emph{Bulletin of The American Mathematical Society} (volume 64:3 (1958)) 
contains surveys of von Neumann's mathematical work, beginning with a valuable article  
 by Ulam~\cite{ulam58}.  The essys by Freeman Dyson~\cite{dyso13} 
 and Paul Halmos~\cite{halm73}  are also excellent.

I thank Caitlin Rizzo, archivist at the Shelby White and Leon Levy Archives Center of the 
Institute for Advanced Study, 
for providing access to its papers on von Neumann.

\end{document}